\documentclass[10pt, reqno]{amsart}
\usepackage{amssymb,latexsym}
\usepackage{eucal}

\newtheorem*{proposition*}{}

\newtheorem*{Lem17.3}{Lemma 17.3}
\newtheorem*{Lem18.3}{Lemma 18.3}
\newtheorem*{Lem18.5}{Lemma 18.5}
\newtheorem*{Lem19.0}{Lemma 19.0}
\newtheorem*{Lem19.1}{Lemma 19.1}
\newtheorem*{Lem19.1.3}{Lemma 19.1.3}

\newtheorem*{Thm}{Theorem}
\newtheorem*{Cor}{Corollary}

\begin{document}

\title[On the asphericity  of  LOT-presentations of groups]
{On the asphericity  of  LOT-presentations of groups}
\author{S.V. Ivanov}
\address{Department of Mathematics\\
University of Illinois \\
Urbana\\   IL 61801\\ U.S.A.} \email{ivanov@math.uiuc.edu}
%Urladdr{}
\thanks{Supported in part by NSF grants  DMS 98-01500, DMS 00-99612}
%\keywords{}
%\date{}

\begin{abstract}
Let $U$ be an arbitrary word in letters $x_1^{\pm 1},  \dots, x_m^{\pm 1}$
and $m  \ge 2$. We prove that the group presentation $\langle   x_1,  \dots,
x_m  \,  \|\ \,  U x_i U^{-1} = x_{i+1}, \, i=1,\dots, m-1  \rangle$ is
aspherical. The proof is based upon prior partial results of A.~Klyachko and
the author on the asphericity of such presentations.
\end{abstract}
\maketitle

Suppose that  $U$ is a word in letters $x_1^{\pm 1},  \dots, x_m^{\pm 1}$,
$m \ge  2$, and a group $G$ is given by a presentation  of the form
\begin{equation}
G = \langle  x_1,  \dots, x_m  \  \|\  \
U x_i U^{-1} = x_{i+1}, \ i=1, \dots, m-1  \rangle .
\end{equation}

It is proved by Klyachko and the author \cite{IK01} that the following claims
(C1)--(C2) hold. Note that for $m=2$ these claims are immediate from classical
Magnus'  and Lyndon's results on one-relator groups, see \cite{MKS66}, \cite{LS77}.

\begin{enumerate}
\item[(C1)]  The presentation (1) of  $G$  is aspherical if
the conjugating word $U$ does not have the form $U_2 U_1$, where $U_1$ is a word
in letters $x_1^{\pm 1},  \dots, x_{m-1}^{\pm 1}$ and $U_2$ is a word in letters
$x_2^{\pm 1},  \dots, x_{m}^{\pm 1}$.

\item[(C2)]   The (images of) letters
$x_1,  \dots, x_{m-1}$  freely generate a free subgroup of the group
$G$ given by presentation (1)
if and only if  the conjugating  word $U$  does not have the form $U_2 U_1$,
where $U_1$ is a word in letters $x_1^{\pm 1},  \dots, x_{m-1}^{\pm 1}$
and $U_2$ is a word in letters $x_2^{\pm 1},  \dots, x_{m}^{\pm 1}$.
\end{enumerate}

The first claim (C1) is of interest because  of
still unsettled Whitehead asphericity conjecture \cite{W41}
that states that every subcomplex of an aspherical 2-complex
is aspherical (see \cite{H83}, \cite{L96}, \cite{I98}, \cite{I99a}, \cite{I99b}),
or, equivalently, every subpresentation
of an aspherical group presentation is aspherical. Observe
that presentation (1) is a subpresentation of
a balanced presentation of the trivial group (which is aspherical)
obtained from (1) by adding a letter to the relator set.
Claim (C1) is also of interest because
the asphericity of presentation (1) is a special case of
the separately conjectured asphericity of LOT-presentations  of groups which,
as was proved by Howie \cite{H85}, is equivalent to the open problem on
the asphericity of ribbon disk complements.

In this note, we will apply a "stabilization" trick to strengthen
claim (C1) of \cite{IK01} by lifting the restriction $U \neq U_2
U_1$.
\begin{Thm}
The group presentation $(1)$ is aspherical for an arbitrary word $U$.
\end{Thm}
\begin{Cor}
Let $k_1, \dots, k_{m-1}$ be some integers,
$U$ be a word in letters $x_1^{\pm 1},  \dots$, $x_m^{\pm 1}$
and $m \ge  2$. Then a presentation  of the form
$
G = \langle  x_1,  \dots, x_m  \  \|\  \ U^{k_i} x_i U^{-k_i} =
x_{i+1}, \ i=1,\dots, m-1  \rangle
$
is apsherical.
\end{Cor}

{\em Proof of Theorem.} Note that for $m=2$ we have a one-relator presentation and
our Theorem follows from classical Lyndon's results on one-relator groups, see \cite{LS77}.
Hence, arguing by induction on $m \ge 2$, in view of claim (C1), we can assume that
$m > 2$ and the conjugating word $U$  has the form $U_2 U_1$, where $U_1$ is a word in letters
$x_1^{\pm 1},  \dots, x_{m-1}^{\pm 1}$ and $U_2$ is a word in letters
$x_2^{\pm 1},  \dots, x_{m}^{\pm 1}$.

Let $S$ be a word in letters $x_1^{\pm 1},  \dots, x_{m-1}^{\pm 1}$. By $S^\alpha$ denote
the word  obtained from $S$ by increasing the index of each letter of $S$ by 1.
Suppose that $X, Y$ are some words  in letters
$x_1^{\pm 1},  \dots, x_{m}^{\pm 1}$
We will write $X \overset G = Y$ if the natural images of
words $X, Y$ in the group $G$ given by (1) are equal.

It easily follows from defining relations of  the group $G$ that
$U U_1 U^{-1} \overset G  = U_1^\alpha $ or
$(U_2 U_1) U_1 (U_2 U_1) ^{-1} \overset G  = U_1^\alpha $, whence
\begin{equation}  %2
U_2 U_1  \overset G  =  U_1^\alpha U_2 .
\end{equation}

Consider another group  $H$  given by the following presentation
\begin{equation}  %3
H = \langle  x_1,  \dots, x_m  \  \|\  \  U_1^\alpha U_2 x_i (U_1^\alpha U_2)^{-1}
= x_{i+1}, \ i=1,\dots, m-1  \rangle .
\end{equation}

Note that
$U_1^\alpha U_2 U_1 (U_1^\alpha U_2)^{-1} \overset H =  U_1^\alpha$ .
Hence,
\begin{equation}  %4
U_2 U_1  \overset H  =  U_1^\alpha  U_2 .
\end{equation}

Consequently, it follows from equalities (2) and (4) that all
defining relations of $H$ hold in $G$ and vice versa, that is,
the groups $G$ and $H$ given by (1) and (3) are naturally isomorphic.

Denote $I_{m-1} = \{1, \dots, m-1\}$ and   consider a sequence of elementary
Andrews-Curtis transformations (AC-moves) applied to presentation (1).
\begin{gather*}
 \langle  x_1,  \dots, x_m  \  \|  \  U_2 U_1 x_i (U_2 U_1)^{-1}
= x_{i+1}, \ i \in I_{m-1}  \rangle  \  \to \\
 \langle  x_1,  \dots, x_m, z  \  \| \    z=1, \ U_2 U_1 x_i (U_2 U_1)^{-1}
= x_{i+1}, \  i \in I_{m-1}  \rangle  \  \to \\
 \langle  x_1,  \dots, x_m, z  \  \| \   z U_1^\alpha U_2(U_2 U_1)^{-1}=1, \
U_2 U_1 x_i (U_2 U_1)^{-1} = x_{i+1}, \ i \in I_{m-1}  \rangle  \  \to \\
K = \langle  x_1,  \dots, x_m, z  \  \| \    z U_1^\alpha U_2(U_2 U_1)^{-1}=1, \
z U_1^\alpha U_2 x_i (z U_1^\alpha U_2)^{-1} = x_{i+1}, \ i \in I_{m-1}  \rangle .
\end{gather*}

Note that it follows from  the last $m-1$ relations of the last
presentation that in the group $K$ given by this presentation we have
$$
 z U_1^\alpha U_2  U_1^{-1} (z U_1^\alpha U_2)^{-1} \overset K  =
(U_1^\alpha)^{-1} .
$$
Therefore,  we can continue the chain of elementary  AC-moves,
replacing the relator $R = z U_1^\alpha U_2(U_2 U_1)^{-1}$ by
$$
\left( z U_1^\alpha U_2  U_1^{-1}
(z U_1^\alpha U_2)^{-1} U_1^\alpha \right)^{-1} R =
\left( R   (U_1^\alpha )^{-1} z^{-1}  U_1^\alpha   \right)^{-1} R =
(U_1^\alpha )^{-1} z   U_1^\alpha
$$
or just by $z$.
\begin{gather*}
\langle  x_1,  \dots, x_m, z  \  \|  \  z U_1^\alpha U_2(U_2 U_1)^{-1}=1 , \
z U_1^\alpha U_2 x_i (z U_1^\alpha U_2)^{-1} = x_{i+1}, \
i \in I_{m-1}  \rangle  \ \to \\
\langle  x_1,  \dots, x_m, z  \  \| \  z=1, \
z U_1^\alpha U_2 x_i (z U_1^\alpha U_2)^{-1} = x_{i+1}, \ i \in I_{m-1}  \rangle \  \to \\
\langle  x_1,  \dots, x_m, z  \  \| \ z=1, \
U_1^\alpha U_2 x_i (U_1^\alpha U_2)^{-1} = x_{i+1}, \ i \in I_{m-1} \rangle \  \to \\
H = \langle  x_1,  \dots, x_m \  \| \
U_1^\alpha U_2 x_i (U_1^\alpha U_2)^{-1} = x_{i+1}, \ i \in I_{m-1}  \rangle  .
\end{gather*}

Now we see that presentation (3)  is Andrews-Curtis  equivalent
to presentation  (1). Since elementary AC-moves preserve the asphericity
of a presentation (see also \cite{I99a}, \cite{I98}), the asphericity of (1) would follow
from the asphericity of (3).  Since the word $U_1^\alpha U_2$ has
no occurrences of $x_1^{\pm 1}$, it follows that the asphericity of (3)
is equivalent to the asphericity of presentation
$$
\langle  x_2,  \dots, x_m   \  \| \
U_1^\alpha U_2 x_i (U_1^\alpha U_2)^{-1} = x_{i+1}, \ i=2, \dots, m-1  \rangle  .
$$
Now it remains to refer to the induction hypothesis and Theorem is proved. \qed
\medskip

To prove Corollary, we note that introduction of new letters $x_{i,j}$
and splitting relations $U^{k_i} x_i U^{-k_i} = x_{i+1}$, where,
say, $k_i >0$,  into several relations
$$
U x_i U^{-1} = x_{i,1}, \  \dots , \ U x_{i, k_{i -1}} U^{-1} = x_{i+1} ,
$$
where $x_{i,1},  \dots , x_{i, k_{i -1}}$  are new letters,
result in a presentation whose asphericity is equivalent
to the asphericity of the original presentation and for which
(after obvious reindexing) all the numbers $k_i$ are equal to $\pm 1$.
Applying evident AC-moves, we can eliminate some  of the letters and relations
and turn all $k_i$ into 1. Now we can refer to proven Theorem. \qed

\end{document}